\documentclass[12pt]{article}
\usepackage[letterpaper, left=1.5in, right=1in,top=1in,bottom=1in]{geometry}
\usepackage{tikz}
\usepackage{chngcntr}
\usepackage{authblk}
\counterwithin{figure}{section}
\usetikzlibrary{arrows}
\usepackage{placeins}	
\usepackage{amsfonts}
\usepackage{graphicx}
\usepackage{standalone}
\usepackage{amsmath, amsthm, amssymb}
\usepackage[capitalize, noabbrev]{cleveref}
\newtheorem{thrm}{Theorem}[section]

\newtheorem{lem}[thrm]{Lemma}
\newtheorem{cor}[thrm]{Corollary}

\newenvironment{pf}           {\noindent{\bf Proof:} }%
                                {\null\hfill$\Box$\par\medskip}
                           
\usepackage[nottoc]{tocbibind}

\newcommand{\pink}[1]{{\color{magenta} #1}}

\begin{document}

\title{Orthogonal Colourings of Random Geometric Graphs}
\author[1]{Jeannette Janssen}
\author[2]{Kyle MacKeigan}
\affil[1]{Department of Mathematics, Dalhousie University}
\affil[ ]{Email: \textit {Jeannette.janssen@dal.ca}}
\affil[2]{Department of Mathematics, Mount Saint Vincent University}
\affil[ ]{Email: \textit{Kyle.m.mackeigan@gmail.com}}
\maketitle

\begin{abstract}
In this paper, we study orthogonal colourings of random geometric graphs. Two colourings of a graph are orthogonal if they have the property that when two vertices receive the same colour in one colouring, then those vertices receive distinct colours in the other colouring.
A random geometric graph $RG(n,r)$ is a graph constructed by randomly placing $n$ vertices in the unit square and connecting two vertices with an edge if and only if their distance is less than the threshold $r$. 
We show first that random geometric graphs with $r>n^{-\alpha}$, where $0\leq \alpha \leq\frac{1}{4}$, have an orthogonal colouring using $n^{1-2\alpha}(1+o(1))$ colours with high probability. Then, we show for an infinite number of values of $n$, random geometric graphs with threshold $r<cn^{-\frac{1}{4}}$, $c<1$, have an optimal orthogonal colouring with high probability. We obtain both of these results by constructing orthogonal colourings of the clique grid graph.
\end{abstract}

\section{Introduction}

Two vertex colourings of a graph are \textit{orthogonal} if they have the property that when two vertices receive the same colour in one colouring, then those vertices must receive distinct colours in the other colouring. An \textit{orthogonal colouring} is a pair of orthogonal vertex colourings. The \textit{orthogonal chromatic number} of a graph $G$, denoted $O\chi(G)$, is the minimum number of colours required for an orthogonal colouring.

Orthogonal colourings were first defined in 1985 by Archdeacon, Dinitz, and Harary in the context of edge colourings \cite{archdeacon1985orthogonal}. Later in 1999, Caro and Yuster studied orthogonal colourings in the context of vertex colourings \cite{caro1999orthogonal}. Then in 2013, Ballif studied upper bounds on collections of orthogonal vertex colourings \cite{ballif2013upper}. In previous work \cite{janssen2020orthogonal}, the authors explored orthogonal colourings of Cayley graphs and Cartesian products of graphs. A part of the results presented here were obtained in the second author's doctoral dissertation \cite{mackeigan2021exploration}.

For a graph $G$ with $n$ vertices, in order for two colourings to be orthogonal, at least $\lceil\sqrt{n}\,\rceil$ colours are required. If $O\chi(G)=\lceil\sqrt{n}\,\rceil$, then $G$ is said to have an \textit{optimal orthogonal colouring}. Graphs with optimal orthogonal colourings are of particular interest due to their applications to independent coverings \cite{mackeigan2021independent} and combinatorial game theory \cite{andres2019orthogonal}. In this paper, optimal orthogonal colourings of random geometric graphs are studied.

The \textit{random geometric graph model}, denoted $RG(n,r)$, is defined as follows. In this model, $n$ points are placed in the unit square, $[0,1]^2$, uniformly at random. Two vertices are then connected by an edge if and only if the Euclidean distance between the two vertices is less than $r$. If $G$ is a graph sampled from the random geometric graph model, then this is denoted by $G\sim RG(n,r)$. Random geometric graphs are interesting to study since they can be used to model real world networks  \cite{kozma2015random}.

We will focus on connected graphs. Penrose \cite{penrose1999k} showed that if $\frac{r^2n}{\ln n}\to \infty$ as $n\to \infty$, then $G\sim RG(n,r)$  is connected with high probability. A random geometric graph $RG(n,r)$ with threshold $r=r(n)$ so that  $\frac{r^2n}{\ln n}\to \infty$ as $n\to\infty$ is called a \textit{dense random geometric graph}. McDiarmid \cite{mcdiarmid2003random} showed that if $G\sim RG(n^2,r)$ is a dense random geometric graph, then $\chi(G)\leq \frac{\sqrt{3}}{2}r^2n^2$ with high probability. This was done by constructing a graph homomorphism into the triangular lattice graph. 

In this paper, we study orthogonal colourings of dense random geometric graphs by constructing a graph homomorphism into a structured graph we call the clique grid graph. In Section 2, we define the clique grid graph and establish results on its orthogonal chromatic number. In Section 3, we use the clique grid graph to bound the orthogonal chromatic number of dense random geometric graphs. In Section 3.1 we show that, if $r=n^{-\alpha}$ where $0\leq \alpha \leq\frac{1}{4}$,  a random geometric graph $G\sim RG(n,r)$ has an orthogonal colouring using $n^{1-2\alpha}(1+o(1))$ colours with high probability. In Section 3.2, we show that, if $r\leq cn^{-\frac{1}{4}})$, $c<1$, then a random geometric graph $G\sim RG(n,r)$ has an optimal orthogonal colouring with high probability.

\section{Orthogonal Colourings of the Clique Grid Graph}

We will prove our results about random geometric graphs by mapping such graphs to graphs that have a highly geometric structure. We will refer to such graphs as \emph{clique grid graphs}. In this section, we define clique grid graphs, and give results on their orthogonal chromatic number.

The clique grid graph  will be constructed by taking the strong graph product of the following graph with itself. Let $H(m,d,t)$ be the graph with vertices labelled $v_i^j$ for $0\leq i<m$ and $0\leq j<t$, where two vertices $v_{i_1}^{j_1}$ and $v_{i_2}^{j_2}$ are adjacent if and only if $|i_1-i_2|\leq d$. That is, $ H(m,d,t)$ is the graph obtained by taking $m$ cliques of size $t$, denoted $C_i$, where all of the vertices in $C_{i_1}$ and $C_{i_2}$ are adjacent to one another if and only if $|i_1-i_2|\leq d$. For example, the graph $H(9,1,2)$ is given in Figure \ref{Figure: Graph Bar H(9,1,2)}; it consists of nine cliques of size two, and only adjacent cliques are connected. 

\begin{figure}[h!]
\centering
\begin{tikzpicture}[line cap=round,line join=round,>=triangle 45,x=1.0cm,y=1.0cm]
\draw [line width=2.pt] (0.,0.)-- (1.,0.);
\draw [line width=2.pt] (0.,0.)-- (1.,1.);
\draw [line width=2.pt] (0.,1.)-- (1.,1.);
\draw [line width=2.pt] (0.,1.)-- (1.,0.);
\draw [line width=2.pt] (1.,0.)-- (2.,0.);
\draw [line width=2.pt] (1.,0.)-- (2.,1.);
\draw [line width=2.pt] (1.,1.)-- (2.,1.);
\draw [line width=2.pt] (1.,1.)-- (2.,0.);
\draw [line width=2.pt] (2.,0.)-- (3.,0.);
\draw [line width=2.pt] (2.,0.)-- (3.,1.);
\draw [line width=2.pt] (2.,1.)-- (3.,1.);
\draw [line width=2.pt] (2.,1.)-- (3.,0.);
\draw [line width=2.pt] (0.,0.)-- (0.,1.);
\draw [line width=2.pt] (1.,0.)-- (1.,1.);
\draw [line width=2.pt] (2.,0.)-- (2.,1.);
\draw [line width=2.pt] (3.,0.)-- (3.,1.);
\draw [line width=2.pt] (3.,0.)-- (4.,1.);
\draw [line width=2.pt] (3.,0.)-- (4.,0.);
\draw [line width=2.pt] (3.,1.)-- (4.,1.);
\draw [line width=2.pt] (3.,1.)-- (4.,0.);
\draw [line width=2.pt] (4.,0.)-- (5.,0.);
\draw [line width=2.pt] (4.,0.)-- (5.,1.);
\draw [line width=2.pt] (4.,1.)-- (5.,1.);
\draw [line width=2.pt] (4.,1.)-- (5.,0.);
\draw [line width=2.pt] (5.,0.)-- (6.,0.);
\draw [line width=2.pt] (5.,0.)-- (6.,1.);
\draw [line width=2.pt] (5.,1.)-- (6.,1.);
\draw [line width=2.pt] (5.,1.)-- (6.,0.);
\draw [line width=2.pt] (6.,0.)-- (7.,0.);
\draw [line width=2.pt] (6.,0.)-- (7.,1.);
\draw [line width=2.pt] (6.,1.)-- (7.,1.);
\draw [line width=2.pt] (6.,1.)-- (7.,0.);
\draw [line width=2.pt] (7.,0.)-- (8.,0.);
\draw [line width=2.pt] (7.,0.)-- (8.,1.);
\draw [line width=2.pt] (7.,1.)-- (8.,1.);
\draw [line width=2.pt] (7.,1.)-- (8.,0.);
\draw [line width=2.pt] (8.,0.)-- (8.,1.);
\draw [line width=2.pt] (7.,1.)-- (7.,0.);
\draw [line width=2.pt] (6.,0.)-- (6.,1.);
\draw [line width=2.pt] (5.,1.)-- (5.,0.);
\draw [line width=2.pt] (4.,0.)-- (4.,1.);

\draw (-0.25,0) node[anchor=north west] {$v_0^0$};
\draw (0.75,0) node[anchor=north west] {$v_1^0$};
\draw (1.75,0) node[anchor=north west] {$v_2^0$};
\draw (2.75,0) node[anchor=north west] {$v_3^0$};
\draw (3.75,0) node[anchor=north west] {$v_4^0$};
\draw (4.75,0) node[anchor=north west] {$v_5^0$};
\draw (5.75,0) node[anchor=north west] {$v_6^0$};
\draw (6.75,0) node[anchor=north west] {$v_7^0$};
\draw (7.75,0) node[anchor=north west] {$v_8^0$};
\draw (-0.25,1.75) node[anchor=north west] {$v_0^1$};
\draw (0.75,1.75) node[anchor=north west] {$v_1^1$};
\draw (1.75,1.75) node[anchor=north west] {$v_2^1$};
\draw (2.75,1.75) node[anchor=north west] {$v_3^1$};
\draw (3.75,1.75) node[anchor=north west] {$v_4^1$};
\draw (4.75,1.75) node[anchor=north west] {$v_5^1$};
\draw (5.75,1.75) node[anchor=north west] {$v_6^1$};
\draw (6.75,1.75) node[anchor=north west] {$v_7^1$};
\draw (7.75,1.75) node[anchor=north west] {$v_8^1$};
\begin{scriptsize}
\draw [fill=black] (0.,0.) circle (2.0pt);
\draw [fill=black] (0.,1.) circle (2.5pt);
\draw [fill=black] (1.,0.) circle (2.5pt);
\draw [fill=black] (1.,1.) circle (2.5pt);
\draw [fill=black] (2.,0.) circle (2.5pt);
\draw [fill=black] (2.,1.) circle (2.5pt);
\draw [fill=black] (3.,0.) circle (2.5pt);
\draw [fill=black] (3.,1.) circle (2.5pt);
\draw [fill=black] (4.,0.) circle (2.5pt);
\draw [fill=black] (4.,1.) circle (2.5pt);
\draw [fill=black] (5.,0.) circle (2.5pt);
\draw [fill=black] (5.,1.) circle (2.5pt);
\draw [fill=black] (6.,0.) circle (2.5pt);
\draw [fill=black] (6.,1.) circle (2.5pt);
\draw [fill=black] (7.,1.) circle (2.5pt);
\draw [fill=black] (7.,0.) circle (2.5pt);
\draw [fill=black] (8.,0.) circle (2.5pt);
\draw [fill=black] (8.,1.) circle (2.5pt);
\end{scriptsize}
\end{tikzpicture}
\caption{$H(9,1,2)$}
  \label{Figure: Graph Bar H(9,1,2)}
\end{figure}

Orthogonal colourings of the graph $H(m,d,t)$ are now considered. These will later be extended to orthogonal colourings of the clique grid graph. {It is easy to see that $H(m,d,t)$ has order $mt$ and clique number $t(d+1)$}. Therefore, $t(d+1)$ is a lower bound on the chromatic number, and thus also the orthogonal chromatic number.
The following theorem shows that either $H(m,d,t)$ has an optimal orthogonal colouring, or the orthogonal chromatic number of $H(m,d,t)$ is at most one more than the clique number. 
{
\begin{thrm}\label{Lemma: Colouring Lattice Component}
For all positive integers $m,d,t$,
\begin{equation*}
     O\chi(H(m,d,t))= \left\{
     \begin{array}{lll}
     t(d+1) &\text{if }m\leq t(d+1), & \text{(Case 1)}\\
       t(d+1)+1 &\text{if }t(d+1)<m\leq t(d+1)^2, & \text{(Case 2)}\\
     \lceil\sqrt{mt}\rceil & \text{otherwise.} &\text{(Case 3)} 
     \end{array}
     \right.
 \end{equation*}
 In particular, if we are in Case 3, then $H(m,d,t)$ has an optimal orthogonal colouring.
\end{thrm}
}

\begin{pf}
Suppose first that we are in Case 1, and let $N=t(d+1)$. According to our condition, $mt\leq N^2$. {Consider the following two colourings $c_1,c_2$ of $H(m,d,t)$. For all $0\leq i< m$ and $0\leq j<t$, let
\begin{eqnarray*}
c_1(v_i^j)=&j+it&(\textrm{mod}~N),\text{ and}\\
c_2(v_i^j)=&\left( j+\left\lfloor\frac{i}{d+1}\right\rfloor\right)(\text{mod}~{t})+it&(\textrm{mod}~N).
\end{eqnarray*}
Here and in the rest of the paper we use $x=a (\text{mod}~{b})$ to mean that $a$ is the unique integer $0\leq a<b$ which is equivalent to $x$ modulo $b$.

First, we will show that these two colourings are proper. Suppose that the vertices $v_{i_1}^{j_1}$ and $v_{i_2}^{j_2}$ are adjacent, that is, $|i_1-i_2|\leq d$. 
Since $(i_1,j_1)\not=(i_2,j_2)$ and $j_1,j_2<t$, we have that $j_1+i_1t\not= j_2-i_2t$. Moreover, we have that 
\begin{equation*}\label{eq:c1}
|j_1+i_1t-j_2-i_2t|\leq |i_1-i_2|t+|j_1-j_2|< t(d+1).
\end{equation*}
Since $N= t(d+1)$, it follows that $c_1(v_{i_1}^{j_1})\neq c_1(v_{i_2}^{j_2})$. This proves that $c_1$ is a proper colouring.

Next we show that $c_2$ is proper. We consider two cases. If $i_1=i_2$, then 
\[
\left| c_2(v_{i_1}^{j_1})- c_2(v_{i_2}^{j_2})\right| = (j_1-j_2)(\text{mod~}{t})<t<N.
\]
Since $j_1\not=j_2$, $c_2(v_{i_1}^{j_1})\neq c_2(v_{i_2}^{j_2})$. Suppose then that $i_1\neq i_2$. Let $s_k=\left( j_k+\left\lfloor\frac{i_k}{d+1}\right\rfloor\right)(\text{mod}~t)$ for $k=1,2$ and assume \emph{wlog} that $s_1\geq s_2$ and thus $0\leq s_1-s_2<t$. Then
\[
\left| s_1+i_1t- (s_2+i_2t) \right|\geq |i_1-i_2|t-|s_1-s_2|\geq t-(s_1-s_2)>0,
\]
and similarly,
\[
\left|s_1+i_1t- (s_2+i_2t) \right|\leq |i_1-i_2|t+|s_1-s_2|< dt+t =N,
\]
and thus $c_2(v_{i_1}^{j_1})\neq c_2(v_{i_2}^{j_2})$. 
We now prove that $c_1$ and $c_2$ are orthogonal. 
Since $N=t(d+1)$ we have a repeating pattern of colours on each clique. Specifically, for each $0\leq b<(d+1)$, $c_1$ and $c_2$ both assign colours $a+bt$, $0\leq a<t$ to cliques $C_{b},C_{b+(d+1)},C_{b+2(d+1)},\dots $, where the last clique to receive this set of colours is { clique $C_{b+s(d+1)}$ where $s$ is the largest integer so that $b+s(d+1)< m$.} 

Now consider vertex $v_i^j$ in clique $i=b+c(d+1)$.  { Since $i<m$, we have that $c\leq \lfloor \frac{m}{d+1}\rfloor\leq t$} by the conditions of Case 1. This vertex receives colours 
\[
c_1(v_i^j)=j +bt+ct(d+1)\,(\text{mod}~N)=j+bt \text{
and }c_2(v_i^j)=(j+c)\,(\text{mod}~{t}) +bt.
\]
As argued earlier, $c<t$. Therefore, for all $0\leq b<d+1$, each pair of colours $(j+bt, j+c+bt\,(\text{mod}~{t}))$, $0\leq c<t$ occurs at most once. Colour pairs $(j+b_1t, j+c+b_2t\,(\text{mod}~{t}))$ never occur. Thus, $c_1$ and $c_2$ constitute an orthogonal colouring. {Since $H(m,d,t)$ has clique number and chromatic number equal to $t(d+1)=N$, this orthogonal colouring is optimal.}

{Suppose next that we are in Case 2, and let $N=t(d+1)$. First we show that $O\chi(H(m,d,t)>N$. Suppose an orthogonal colouring $c_1,c_2$ with $N$ colours exists. Then such a colouring must have the repeating pattern as in Case 1.
Specifically, for each $0\leq b<(d+1)$, $c_1$ and $c_2$ must both assign colours $a+bt$, $0\leq a<t$ to cliques $C_{b},C_{b+(d+1)},C_{b+2(d+1)},\dots, C_{b+s(d+1)} $,  where $s$ is the largest integer so that $b+s(d+1)< m$. Namely, cliques $C_0,C_1,\dots ,C_{d}$ form a clique, so admit a unique colouring. Every vertex in $C_{d+1}$ is adjacent to all vertices in clique $C_1,\dots ,C_{d}$, so the only colours available are the colours used on clique $C_0$. Similarly, for each $d+1\leq k<m$, the only colours available to colour clique $C_k$ are the colours used to colour $C_{k-(d+1)}$. Now consider the colourings assigned by  $c_1$ and $c_2$ to cliques $C_0,C_{(d+1)},\dots ,C_{s(d+1)}$, where $s=\lceil 
\frac{m}{d+1}\rceil -1$. Each of these colourings assigns the same set of $t$ colours. There are at most $t^2$ colour pairs with this set of colours, so these cliques together can contain at most $t^2$ vertices, Therefore, we must have that $s<t$ and thus $m\leq t(d+1)$.

To show that $O\chi(H(m,d,t)=N+1$, consider the following two colourings using $N+1$ colours. 
For all $0\leq i< m$ and $0\leq j<t$, let
\begin{eqnarray*}
c_1(v_i^j)=&j+it&(\textrm{mod}~N),\text{ and}\\
c_2(v_i^j)=&( j+it)\,(\text{mod}~{N})+\lfloor\frac{i}{d+1}\rfloor&(\textrm{mod}~N+1).
\end{eqnarray*}

Colouring $c_1$ is the same as in the previous case, and to prove it is proper we did not need the condition on $m$. 
Next we show that $c_2$ is proper. Let $v_{i_1}^{j_1}$ and $ v_{i_2}^{j_2}$ be two adjacent  vertices. We consider two cases. If $i_1=i_2$, then $( j_1+i_1t)(\text{mod}~{N})\not= ( j_2+i_2t)(\text{mod}~{N})$ as argued before, and thus $c_2(v_{i_1}^{j_1})\not = c_2(v_{i_2}^{j_2})$.
Suppose then that $0<i_1-i_2\leq d$. Then
\begin{equation}\label{eq:remainders}
%|c_2(v_{i_1}^{j_1}) -c_2(v_{i_2}^{j_2})|=
\left\lfloor\frac{i_1}{d+1}\right\rfloor-\left\lfloor\frac{i_2}{d+1}\right\rfloor \leq \left\lfloor\frac{i_1-i_2}{d+1}\right\rfloor +1\leq 1.
%\,(\textrm{mod}~N+1)\not= 0.
\end{equation}
If the quantity in \eqref{eq:remainders} equals zero, then $c_2(v_{i_k}^{j_k})=c_1(v_{i_k}^{j_k})$ for $k=1,2$, which is impossible since $c_1$ is proper. Then
\[
1\leq(j_1+i_1t-j_2-i_2t)\,(\text{mod}~{N})+\left\lfloor\frac{i_1}{d+1}\right\rfloor-\left\lfloor\frac{i_2}{d+1}\right\rfloor\leq (N-1)+1\leq N,
\]
and thus $c_2(v_{i_1}^{j_1})\not = c_2(v_{i_2}^{j_2})$.

We now prove that $c_1$ and $c_2$ are orthogonal. Suppose $c_1(v_{i_1}^{j_1})=c_1(v_{i_2}^{j_2})$. Since $c_1$ is a proper colouring, we must have that $|i_1-i_2|\geq d+1$. Since $i_k<m\leq N(d+1)$ for $k=1,2$,
\[
|c_2(v_{i_1}^{j_1}) -c_2(v_{i_2}^{j_2})|=\left|\left\lfloor\frac{i_1}{d+1}\right\rfloor-\left\lfloor\frac{i_2}{d+1}\right\rfloor\right| \,(\textrm{mod}~N+1)\not= 0.
\]
}

Suppose finally that we are in Case 3 and let $N=\lceil\sqrt{mt}\rceil$
. 
By the conditions of this case, $mt>t^2(d+1)^2$, so $N\geq t(d+1)+1$. Consider the following two colourings $c_1,c_2$ of $H(m,d,t)$. For all $0\leq i< m$ and $0\leq j<t$, let
\begin{eqnarray*}
c_1(v_i^j)=&j+it&(\textrm{mod}~N),\text{ and}\\
c_2(v_i^j)=&j+it+\left\lfloor\frac{j+it}{N}\right\rfloor&(\textrm{mod}~N).
\end{eqnarray*}
Note that $j+it<mt$.

{Colouring $c_1$ is the same as in Case 1, so we proceed to show that $c_2$ is proper. Suppose that $v_{i_1}^{j_1}\not= v_{i_2}^{j_2}$ and $0\leq i_1-i_2\leq d$. 

Since $N>t(d+1)$, we have that 
\begin{eqnarray*}
   \left\lfloor\frac{j_1+i_1t}{N}\right\rfloor-\left\lfloor\frac{j_2+i_2t}{N}\right\rfloor
   &\leq & \left\lfloor\frac{ j_1+i_1t-(j_2+i_2t)}{N}\right\rfloor+1\\
   &\leq &\left\lfloor\frac{| j_1-j_2|+(i_1-i_2)t}{N}\right\rfloor+1\\
   &\leq & \left\lfloor\frac{t-1+dt}{N}\right\rfloor+1\leq 1.
\end{eqnarray*}

It follows that
\begin{eqnarray*}
j_1+i_1t+\left\lfloor\frac{j_1+i_1t}{N}\right\rfloor-\left(j_2+i_2t+\left\lfloor\frac{j_2+i_2t}{N}\right\rfloor\right)
&\leq &|j_1-j_2|+|i_1-i_2|t+1
\\
&\leq & (t-1)+dt+1=t(d+1)<N,
\end{eqnarray*}}
and
\begin{eqnarray*}
j_1+i_1t+\left\lfloor\frac{j_1+i_1t}{N}\right\rfloor-\left(j_2+i_2t+\left\lfloor\frac{j_2+i_2t}{N}\right\rfloor\right)
&\geq &(j_1-j_2)+ (i_1-i_2)t+1\\
&\geq & -(t-1)+t+0>0.
\\
\end{eqnarray*}}
 It follows that $c_2(v_{i_1}^{j_1})\neq c_2(v_{i_2}^{j_2})$. Hence, $c_1$ and $c_2$ are proper colourings of $H(m,d,t)$.  

We will now show that $c_1$ and $c_2$  are orthogonal colourings. 
Suppose that $c_1(v_{i_1}^{j_1})=c_1(v_{i_2}^{j_2})$ and $c_2(v_{i_1}^{j_1})=c_2(v_{i_2}^{j_2})$ where $i_1\neq i_2$ or $j_1\neq j_2$. Since $c_1(v_{i_1}^{j_1})=c_1(v_{i_2}^{j_2})$, this implies that $i_1t+j_1=i_2t+j_2+cN$. {Since $0\leq j_k<t$ and $0\leq i_k<m$ for $k=1,2$ we have that $0<(j_1-j_2)+(i_1-i_2)t<mt\leq N^2$, so $0<c<N$.} It follows that
\begin{align*}
c_2(v_{i_1}^{j_1})&=i_2t+j_2+\left\lfloor\frac{i_2t+j_2+cN}{N}\right\rfloor(\textrm{mod}~N)\\
&=i_2t+j_2+c+\left\lfloor\frac{i_2t+j_2}{N}\right\rfloor(\textrm{mod}~N)\\
&=c_2(v_{i_2}^{j_2})+c\, (\textrm{mod}~N).
\end{align*}
Since $c_2(v_{i_1}^{j_1})=c_2(v_{i_2}^{j_2})$, this gives that $c\equiv 0(\textrm{mod}~N)$, contradicting {$0<c<N$}.

\end{pf}

Now, we can define the clique grid graph and extend the orthogonal colouring of $H(m,d,t)$. The clique grid graph, denoted $L(m^2,d,t^2)$, is obtained by taking the strong product of $H(m,d,t)$ with itself. That is, $L(m^2,d,t^2)= H(m,d,t)\boxtimes H(m,d,t)$. Alternatively, $L(m^2,d,t^2)$ can be viewed as $m^2$ cliques of size $t^2$, denoted $C_{i,j}$, where all of the vertices in $C_{i_1,j_1}$ and $C_{i_2,j_2}$ are adjacent to one another if and only if $|i_1-i_2|\leq d$ and $|j_1-j_2|\leq d$. For example, the clique grid $L(25,1,1)$ 
%\pink{I would prefer a picture where $t>1$, for example $L(9, 1, 4)$, since the case $t=1$ is just the strong product of paths, so does not show the particular clique structure.}
is given in Figure \ref{Figure: Clique Grid}. 

\begin{figure}[h!]
\centering
\begin{tikzpicture}[line cap=round,line join=round,>=triangle 45,x=1.0cm,y=1.0cm]
\draw [line width=2.pt] (0.,4.)-- (1.,4.);
\draw [line width=2.pt] (1.,4.)-- (2.,4.);
\draw [line width=2.pt] (2.,4.)-- (3.,4.);
\draw [line width=2.pt] (3.,4.)-- (4.,4.);
\draw [line width=2.pt] (4.,3.)-- (3.,3.);
\draw [line width=2.pt] (3.,3.)-- (2.,3.);
\draw [line width=2.pt] (2.,3.)-- (1.,3.);
\draw [line width=2.pt] (1.,3.)-- (0.,3.);
\draw [line width=2.pt] (0.,2.)-- (1.,2.);
\draw [line width=2.pt] (1.,2.)-- (2.,2.);
\draw [line width=2.pt] (2.,2.)-- (3.,2.);
\draw [line width=2.pt] (3.,2.)-- (4.,2.);
\draw [line width=2.pt] (4.,1.)-- (3.,1.);
\draw [line width=2.pt] (3.,1.)-- (2.,1.);
\draw [line width=2.pt] (2.,1.)-- (1.,1.);
\draw [line width=2.pt] (1.,1.)-- (0.,1.);
\draw [line width=2.pt] (0.,0.)-- (1.,0.);
\draw [line width=2.pt] (1.,0.)-- (2.,0.);
\draw [line width=2.pt] (2.,0.)-- (3.,0.);
\draw [line width=2.pt] (3.,0.)-- (4.,0.);
\draw [line width=2.pt] (0.,4.)-- (0.,3.);
\draw [line width=2.pt] (0.,3.)-- (0.,2.);
\draw [line width=2.pt] (0.,2.)-- (0.,1.);
\draw [line width=2.pt] (0.,1.)-- (0.,0.);
\draw [line width=2.pt] (1.,4.)-- (1.,3.);
\draw [line width=2.pt] (1.,3.)-- (1.,2.);
\draw [line width=2.pt] (1.,2.)-- (1.,1.);
\draw [line width=2.pt] (1.,1.)-- (1.,0.);
\draw [line width=2.pt] (2.,0.)-- (2.,1.);
\draw [line width=2.pt] (2.,1.)-- (2.,2.);
\draw [line width=2.pt] (2.,2.)-- (2.,3.);
\draw [line width=2.pt] (2.,3.)-- (2.,4.);
\draw [line width=2.pt] (3.,4.)-- (3.,3.);
\draw [line width=2.pt] (3.,3.)-- (3.,2.);
\draw [line width=2.pt] (3.,2.)-- (3.,1.);
\draw [line width=2.pt] (3.,1.)-- (3.,0.);
\draw [line width=2.pt] (4.,0.)-- (4.,1.);
\draw [line width=2.pt] (4.,1.)-- (4.,2.);
\draw [line width=2.pt] (4.,2.)-- (4.,3.);
\draw [line width=2.pt] (4.,4.)-- (4.,3.);
\draw [line width=2.pt] (0.,4.)-- (1.,3.);
\draw [line width=2.pt] (0.,3.)-- (1.,4.);
\draw [line width=2.pt] (1.,4.)-- (2.,3.);
\draw [line width=2.pt] (1.,3.)-- (2.,4.);
\draw [line width=2.pt] (2.,4.)-- (3.,3.);
\draw [line width=2.pt] (2.,3.)-- (3.,4.);
\draw [line width=2.pt] (3.,4.)-- (4.,3.);
\draw [line width=2.pt] (3.,3.)-- (4.,4.);
\draw [line width=2.pt] (0.,3.)-- (1.,2.);
\draw [line width=2.pt] (0.,2.)-- (1.,3.);
\draw [line width=2.pt] (1.,3.)-- (2.,2.);
\draw [line width=2.pt] (1.,2.)-- (2.,3.);
\draw [line width=2.pt] (2.,3.)-- (3.,2.);
\draw [line width=2.pt] (2.,2.)-- (3.,3.);
\draw [line width=2.pt] (3.,3.)-- (4.,2.);
\draw [line width=2.pt] (3.,2.)-- (4.,3.);
\draw [line width=2.pt] (0.,2.)-- (1.,1.);
\draw [line width=2.pt] (0.,1.)-- (1.,2.);
\draw [line width=2.pt] (1.,2.)-- (2.,1.);
\draw [line width=2.pt] (1.,1.)-- (2.,2.);
\draw [line width=2.pt] (2.,2.)-- (3.,1.);
\draw [line width=2.pt] (2.,1.)-- (3.,2.);
\draw [line width=2.pt] (3.,2.)-- (4.,1.);
\draw [line width=2.pt] (3.,1.)-- (4.,2.);
\draw [line width=2.pt] (0.,1.)-- (1.,0.);
\draw [line width=2.pt] (0.,0.)-- (1.,1.);
\draw [line width=2.pt] (1.,1.)-- (2.,0.);
\draw [line width=2.pt] (1.,0.)-- (2.,1.);
\draw [line width=2.pt] (2.,1.)-- (3.,0.);
\draw [line width=2.pt] (2.,0.)-- (3.,1.);
\draw [line width=2.pt] (3.,1.)-- (4.,0.);
\draw [line width=2.pt] (3.,0.)-- (4.,1.);
\begin{scriptsize}
\draw [fill=black] (0.,0.) circle (2.0pt);
\draw [fill=black] (1.,0.) circle (2.5pt);
\draw [fill=black] (0.,1.) circle (2.5pt);
\draw [fill=black] (0.,2.) circle (2.5pt);
\draw [fill=black] (0.,3.) circle (2.5pt);
\draw [fill=black] (0.,4.) circle (2.5pt);
\draw [fill=black] (1.,4.) circle (2.5pt);
\draw [fill=black] (2.,4.) circle (2.5pt);
\draw [fill=black] (3.,4.) circle (2.5pt);
\draw [fill=black] (4.,4.) circle (2.5pt);
\draw [fill=black] (4.,3.) circle (2.5pt);
\draw [fill=black] (3.,3.) circle (2.5pt);
\draw [fill=black] (2.,3.) circle (2.5pt);
\draw [fill=black] (1.,3.) circle (2.5pt);
\draw [fill=black] (1.,2.) circle (2.5pt);
\draw [fill=black] (2.,2.) circle (2.5pt);
\draw [fill=black] (3.,2.) circle (2.5pt);
\draw [fill=black] (4.,2.) circle (2.5pt);
\draw [fill=black] (1.,1.) circle (2.5pt);
\draw [fill=black] (2.,1.) circle (2.5pt);
\draw [fill=black] (3.,1.) circle (2.5pt);
\draw [fill=black] (4.,1.) circle (2.5pt);
\draw [fill=black] (4.,0.) circle (2.5pt);
\draw [fill=black] (3.,0.) circle (2.5pt);
\draw [fill=black] (2.,0.) circle (2.5pt);
\end{scriptsize}
\end{tikzpicture}
\caption{$L(25,1,1)$}
\label{Figure: Clique Grid}
\end{figure}

By using the orthogonal colourings of $H(m,d,t)$ from Lemma \ref{Lemma: Colouring Lattice Component}, an orthogonal colouring of $L(m^2,d,t^2)$ is obtained. This comes from the fact that for any graphs $G$ and $H$, {we can use a previous result, proved in \cite{mackeigan2020orthogonal}
along with other graph product results. 
\begin{thrm}[\cite{mackeigan2020orthogonal}] For any two graphs $G$ and $H$,
\[
O\chi(G\boxtimes H)\leq O\chi(G)O\chi(H).
\]
\end{thrm}

Therefore, we have  the following corollary.

\begin{cor}\label{Lemma: Orthogonal Colouring Lattice}
For all positive integers $m,d,t$,
\begin{equation*}
     O\chi(L(m^2,d,t^2))\leq \left\{
     \begin{array}{lll}
     t^2(d+1)^2 &\text{if }m\leq t(d+1), & \text{(Case 1)}\\
          (t(d+1)+1)^2 &\text{if } t(d+1)<m\leq t(d+1)^2, & \text{(Case 2)}\\
     \left(\lceil\sqrt{mt}\rceil\right)^2 & \text{otherwise.} &\text{(Case 3)} 
     \end{array}
     \right.
 \end{equation*}
    \end{cor}}

In particular, if we are in Case 3 and $mt$ is a square, then $L(m^2,d,t^2)$ has an optimal orthogonal colouring. Note also that $\omega (L^2(m^2,d,t^2))=t^2(d+1)^2$, and thus we have equality in Case 1. 

In the next section, we show that an injective graph homomorphism from dense random geometric graphs into the clique grid graph exists with high probability. If such a graph homomorphisms exists, then the orthogonal colourings of $L(m^2,d,t^2)$ in Lemma \ref{Lemma: Orthogonal Colouring Lattice} can be applied to the dense random geometric graphs.

\section{Dense Random Geometric Graphs}

To obtain orthogonal colourings of dense random geometric graphs, we will show that with high probability and for the appropriate choice of parameters, $G\sim RG(n,r)$ is isomorphic to a subgraph of $L(m^2,d,t^2)$. We distinguish two cases. First we consider $RG(n,r)$ where $r=n^{-\alpha}$, $0<\alpha\leq 1/4$. We will see that the orthogonal chromatic number in this case is mostly constrained by the clique number. After that we will see that if $\alpha>1/4$, then the graph has, with high probability, an optimal orthogonal colouring.

\subsection{Orthogonal colourings close to the clique number}
{To obtain a subgraph isomorphism between a random geometric graph and a clique grid graph, the unit square is divided into $m\times m$ equal size squares. In particular, for $l=\frac{1}{m}$, the set $S_{ij}$ will contain the vertices of $G$ in the square with coordinates in the subsquare $((i-1)l,il)\times ((j-1)l,jl)$.

To show that $G$ is isomorphic to a subgraph of $L(m^2,d,t^2)$, all of the vertices in $S_{ij}$ are mapped to vertices in the cliques $C_{ij}$ in $L(m^2,d,t^2)$.} To show that this is a subgraph isomorphism with $H\subseteq L(m^2,d,t^2)$, we show that with high probability, for all $i,j$, that $|S_{ij}|\leq |C_{ij}|=t^2$. Additionally, we show that if two vertices are adjacent in $G$, then their images in $H$ are adjacent.

Now, notice that two vertices $u\in C_{i_1,j_1}$ and $v\in C_{i_2,j_2}$ are adjacent if and only if $|i_1-i_2|\leq d$ and $|j_1-j_2|\leq d$. On the other hand, two vertices $u\in S_{i_1,j_1}$ and $v\in S_{i_2,j_2}$ are adjacent in $G$ if and only if their Euclidean distance is less than $r$. To distinguish between Euclidean distance and the absolute value, the \textit{Euclidean distance} between two points $u$ and $v$ is denoted by $\| u-v\|$.

We will show in Lemma \ref{Lemma: Edge Correctness} that $\|u-v\pink{}|<r$ implies that $|i_1-i_2|<\frac{r}{l}+1$ and $|j_1-j_2|<\frac{r}{l}+1$. Therefore, we define $r=n^{-\alpha},l=\left\lceil\frac{\sqrt{n}}{\ln n}\right\rceil^{-1} $, and $d=\left\lceil\frac{n^{\frac{1}{2}-\alpha}}{\ln n}\right\rceil+2 $ so that
\begin{align}\label{eq:r/l+1<d}
\frac{r}{l}+1&=n^{-\alpha}\left\lceil\frac{\sqrt{n}}{\ln n}\,\right\rceil+1\notag\\
&\leq \frac{n^{1/2-\alpha}}{\ln n}+n^{-\alpha}+1\notag\\
&< \left\lceil\frac{n^{1/2-\alpha}}{\ln n}\,\right\rceil+2= d
\end{align}

Therefore, Lemma \ref{Lemma: Edge Correctness} will give that if $\|u-v\|<r$ then $|i_1-i_2|<\frac{r}{l}+1<d$. Hence, the subgraph isomorphism described will preserve the edges. Additionally, we define two other parameters, $t=\lceil\ln n\rceil$ and $m=\left\lceil\frac{\sqrt{n}}{\ln n}\right\rceil$. These two parameters are defined in this way to satisfy the probabilistic lemmas proved later. For reference, the follow parameters are used throughout this section.

\begin{align}
t&=\lceil\ln n\rceil  \label{Equation: parameter t}\\
m&=\left\lceil\frac{\sqrt{n}}{\ln n}\right\rceil  \label{Equation: parameter m}\\
l&=\left\lceil\frac{\sqrt{n}}{\ln n}\right\rceil^{-1}  \label{Equation: parameter l}\\
d&=\left\lceil\frac{n^{\frac{1}{2}-\alpha}}{\ln n}\right\rceil+2  \label{Equation: parameter d}
%r&=n^{-\alpha}  \label{Equation: parameter r}
\end{align}

Recall that an event $E$ occurs \textit{with high probability} if as $n$ tends to infinity, the probability that $E$ occurs tends to one. We will show that with high probability and for all $i,j$, that $|S_{ij}|\leq |C_{ij}|=t^2$. We prove this result with the well-known Chernoff's bound, which is now stated.

%\pink{Reference for Chernoff bound not working}
\begin{thrm}[Chernoff's Bound] %\cite{chernoff2014career}]
\label{ChernoffBound}
Suppose that $X_1,X_2,\dots ,X_n$ are independent random variables taking values in $\{0,1\}$. Let $X$ denote their sum and let $\mu=\mathbb{E}[X]$. For any $\delta\geq 0$, it follows that $$\mathbb{P}(X<(1-\delta)\mu)\leq e^{-\frac{\delta^2\mu}{2}}\text{ and }\mathbb{P}(X>(1+\delta)\mu)\leq e^{-\frac{\delta^2\mu}{3}}$$
\end{thrm}

Chernoff's bound gives exponentially decreasing bounds and can be applied to bound the probability that $|S_{ij}|\leq |C_{ij}|$ for a single square. However, it is required that this inequality holds for all squares in the partition. To extend this result to all squares, the Union Bound is required. For the events $E_1,\dots, E_n$, let $\cup_i E_i$ denote the event that at least one of the events occurs. The Union Bound is then as follows:
$$\mathbb{P}(\cup_i E_i)\leq \sum_i^n \mathbb{P}(E_i).$$

\begin{lem}\label{Lemma: Chernoff Partition}
For $G\sim RG(n,r)$. Let $t,m,l$ be the parameters in Equations \ref{Equation: parameter t}, \ref{Equation: parameter m} and \ref{Equation: parameter l}.  For all $1\leq i,j\leq m$, let $S_{ij}$ denote the set of vertices of $G$ in the square $((i-1)l,il)\times ((j-1)l,jl)$. With high probability and for all $i,j$, $|S_{i,j}|\leq t^2$.
\end{lem}

\begin{pf}
First, fix the indices $i$ and $j$. For all vertices $v\in V(G)$, define the random variable $X_v$ as $X_v=0$ if $v\not\in S_{i,j}$ and $X_v=1$ if $v\in S_{i,j}$. Let $X_{i,j}$ denote the sum of the random variables. That is, 
$$
X_{i,j}=\sum_{v\in V} X_v=|S_{i,j}|$$

Recall that the points of $G$ are placed uniformly at random in the unit square, which has an area of 1. Also, the area of each of the $m^2$ subsquares in the partition of the unit square is $1/m^2=l^2=\left\lceil\frac{\sqrt{n}}{\ln n}\right\rceil^{-2}$. Therefore, 
$$\mathbb{E}(X_v)=\mathbb{P}(X_v=1)=\left\lceil\frac{\sqrt{n}}{\ln n}\right\rceil^{-2}.
$$

Since there are $n$ points, the expected number of points in the fixed square is given by $\mu=\mathbb{E}(X_{ij})=n\left\lceil\frac{\sqrt{n}}{\ln n}\right\rceil^{-2}\leq (\ln n)^2=t^2$. Lastly, let $\delta=\frac{\sqrt{3\ln n}}{\sqrt{n}}\left\lceil\frac{\sqrt{n}}{\ln n}\right\rceil>0$. By applying Chernoff's bound, it follows that
\begin{align*}
\mathbb{P}(X_{i,j}>(1+\delta)\mu)&\leq e^{-\frac{\left(\frac{\sqrt{3\ln n}}{\sqrt{n}}\left\lceil\frac{\sqrt{n}}{\ln n}\right\rceil\right)^2 n\left\lceil\frac{\sqrt{n}}{\ln n}\right\rceil^{-2}}{3}}= e^{-\ln n}=n^{-1}
\end{align*}

Thus with high probability and for a fixed $i,j$, $|S_{i,j}|=\mathbb{E}(X_{ij})\leq t^2$. However, to obtain our result, it is required that with high probability and for all $i,j$, that $|S_{i,j}|\leq t^2$. We obtain this by applying the Union Bound:
\begin{align*}
\mathbb{P}(\cup_{i,j}\{X_{i,j}>(1+\delta)\mu\})&\leq \sum_{i=1}^m\sum_{j=1}^m\mathbb{P}(X_{i,j}>(1+\delta)\mu)\\
&\leq m^2n^{-1}\\
&= \left\lceil\frac{\sqrt{n}}{\ln n}\right\rceil^2n^{-1}\\
&\leq(\ln n)^{-2}+2(\sqrt{n}\ln n)^{-1}+n^{-1}
\end{align*} 

Note that each of terms in this expression tends to zero as $n$ tends to infinity. Therefore with high probability and for all $i,j$, $|S_{i,j}|\leq t^2$.
\end{pf}

Lemma \ref{Lemma: Chernoff Partition} gives that the vertices in $S_{i,j}$ in the random geometric graph can be mapped injectively into the cliques $C_{i,j}$ of $L(m^2,d,t^2)$. Thus, one required property of the subgraph isomorphism is obtained. It remains to show that the edges are preserved under this mapping. The following lemma provides this property.

\begin{lem}\label{Lemma: Edge Correctness}
Let $G\sim RG(n,r)$ with $r=n^{-\alpha}$. Let $t,m,l,d$ be the parameters in Equations \ref{Equation: parameter t}, \ref{Equation: parameter m}, \ref{Equation: parameter l}, and \ref{Equation: parameter d}.  For all $1\leq i,j\leq m$, let $S_{ij}$ denote the vertices of $G$ in the square $((i-1)l,il)\times ((j-1)l,jl)$. For $u,v\in V(G)$, suppose that $u\in S_{i_1,j_1}$ and $v\in S_{i_2,j_2}$. If $uv\in E(G)$, then $|i_1-i_2|<d$ and $|j_1-j_2|<d$.
\end{lem}
\begin{pf}
Recall that $d$ is the parameter that provides which cliques in $L(m^2,d,t^2)$ are adjacent. Suppose that $u\in S_{i_1,j_1}$, $v\in S_{i_2,j_2}$, and $uv\in E(G)$. Note that, for any $j$ because there are at least $|i_1-i_2|-1$ inclusive squares between $S_{i_1,j}$ and $S_{i_2,j}$, each of width $l$. Therefore, $|i_1-i_2|<\frac{r}{l}+1$.   Similarly, for any $i$, there are at least $|j_1-j_2|-1$ squares between $S_{i,j_1}$ and $S_{i,j_2}$, each of height $l$. Therefore, $|j_1-j_2|<\frac{r}{l}+1$. However, by the specific choice of parameters, $\frac{r}{l}+1< d$, as in\eqref{eq:r/l+1<d}. Therefore, $|i_1-i_2|<d$ and $|j_1-j_2|<d$.
\end{pf}

Lemma \ref{Lemma: Edge Correctness} provides that the previously described mapping will preserve all of the edges. Therefore, by combining Lemma \ref{Lemma: Chernoff Partition} and Lemma \ref{Lemma: Edge Correctness} with Corollary \ref{Lemma: Orthogonal Colouring Lattice}, we obtain the following result.

\begin{thrm}\label{Theorem: Random Main Result}
Let $t,m,l,d$ be the parameters in Equations \ref{Equation: parameter t}, \ref{Equation: parameter m}, \ref{Equation: parameter l}, and \ref{Equation: parameter d}. For $G\sim RG(n,r)$ where $r=n^{-\alpha}$, if $0\leq\alpha\leq \frac{1}{4}$, then with high probability $$\frac{\sqrt{3}}{2}n^{1-2\alpha}\leq O\chi(G)\leq n^{1-2\alpha}(1+o(1))$$ In particular, if $\alpha=\frac{1}{4}$, then with high probability $$O\chi(G)= \sqrt{n}(1+o(1))$$
\end{thrm}
\begin{pf}
Consider partitioning the unit square into $m\times m$ equal size squares. Let $S_{ij}$ denote the vertices of $G$ in the square with dimensions $((i-1)l,il)\times ((j-1)l,jl)$. Consider the mapping that takes all of the vertices in $S_{i,j}$ and maps them to vertices in $C_{i,j}$ in $L(m^2,d,t^2)$. Lemma \ref{Lemma: Chernoff Partition} and Lemma \ref{Lemma: Edge Correctness} then give that, with high probability, $G$ is a subgraph of $L(m^2,d,t^2)$ through this mapping. To apply Corollary \ref{Lemma: Orthogonal Colouring Lattice} (Case 1 or 2) to find an orthogonal colouring of $L(m^2,d,t^2)$, it is required that {$m\le t(d+1)^2$}.
%$d\geq \frac{\sqrt{m}}{\sqrt{t}}$. 
By the choice of parameters, it follows that

\begin{align*}
(d+1)^2 =\left(\frac{n^{1/2-\alpha}}{\ln n}+3\right)^2&\geq 
\left(\displaystyle{\frac{n^{1/4}}{\ln n}}+3\right)^2
&\text{(since $0\leq \alpha\leq \frac{1}{4}$)}\\
&\geq \displaystyle{\frac{n^{1/2}}{(\ln n)^2}}+9  \\
&= \displaystyle{\frac{\frac{\sqrt{n}}{\ln n}}{\ln n}}+9\\
&\geq \displaystyle{\frac{m}{t}}.
\end{align*}

Therefore, since $m\leq t(d+1)^2$, we are in Case 1 or 2, and Corollary \ref{Lemma: Orthogonal Colouring Lattice} can be applied to find an orthogonal colouring of $L(m^2,d,t^2)$. By substituting the parameters into the bound from Corollary \ref{Lemma: Orthogonal Colouring Lattice}, it follows that
\begin{align*}
O\chi(L(m^2,d,t^2))&\leq (t(d+1)+1)^2\\
&\leq \left((\ln n+1)\left(\frac{n^{\frac{1}{2}-\alpha}}{\ln n}+2\right)+1\right)^2\\
&\leq n^{1-2\alpha}(1+o(1))
\end{align*}

% \pink{jj: adjusted to new bound; original still there, commented out.}
% \begin{align*}
% O\chi(G)\leq O\chi(L(m^2,d,t^2))&= t(d+1)^2\\
% &\leq (\ln n+1)\left(\frac{n^{\frac{1}{2}-\alpha}}{\ln n}+4\right)^2\\
%  &= n^{1-2\alpha}(1+o(1))
%  \end{align*}

When $\alpha=1/4$, the upper bound approximately matches the lower bound of $\sqrt{n}$, so this gives that $O\chi(G)\leq \sqrt{n}(1+o(1))$. On the other hand,
it is known \cite{mcdiarmid2003random} that with high probability, $\chi(G)=\frac{\sqrt{3}}{2}n^{1-2\alpha}$. Therefore, since the orthogonal chromatic number is at least the chromatic number, for $\alpha<1/4$ we have that $ O\chi(G)=c n^{1-2\alpha}(1+o(1))$ with high probability, where $\frac{\sqrt{3}}{2}\leq c\leq 1$.
\end{pf}

\subsection{Optimal orthogonal colourings when $r<n^{-1/4}$}

We now consider the case where $r\leq cn^{-1/4}$, $c\in (0,1)$, and establish that in this case we have an optimal orthogonal colouring. Again, we proceed by partitioning the unit square into a similar structure as the clique grid. Whereas in the previous case we had an approximately equal number of points in each cell, in this case we wish to control the number of points in each cell exactly. The trade-off is that the cells will not be of equal size, and may be rectangular rather than square.  

We assume without loss of generality that $n$ is a perfect square. Since we are considering the dense case, we will assume that $n$ is large enough so that $r\leq (\ln{n})^2 n^{-1/2}$. Let $G\sim RG(n,r)$ and let $P=V(G)$. So $P$ is a collection of points in the unit square, selected independently and uniformly.
Let $n=m^2t^2$ for some integers {$1\leq t\leq m$ (we will later add restrictions on $t$ and $m$).}
Divide the unit square into $m^2$ cells containing exactly $t^2$ points from $P$ in the following way. 
%\pink{(jj: I changed the parameters to $m$ and $t$ to be consistent with the definition of $L(m^2,d,t^2)$})
First, order all of the $n$ points in $P$ according to their $y$-values. That is, suppose that $y_1<y_2<\dots<y_{n}$. Then for $0\leq i\leq m$,  define $y_i^*$  as

\begin{equation}\label{def: yi*} 
y_i^*=\left\{
\begin{array}{ll}
    0 & \text{ if } i=0\\
      {y_{itn}} & \text{ if } 1\leq i<{m} \\
      1 & \text{ if } i=m
\end{array} 
\right. 
\end{equation}
 
Now, let $S_i=\{(x,y)\in [0,1]^2:y_{i-1}^*<y\leq y_i^*\}$. By definition, the sets $S_i$ partition the unit square into $m$ strips, each containing $tn=m t^2$ vertices. It is important to note that the $x$-values of the points in each strip are still random. This is because they are independent of the $y$-values. Now, for each $i, 1\leq i\leq t$, we will order the $x$-values of the points in each of the sets $S_i$. That is, we suppose that $x_{i,1}<x_{i,2}<\dots<x_{i,mt^2}$. Then, define $x_{i,j}^*$, $0\leq j\leq m$, as

\begin{equation}\label{def:xij*}
x_{i,j}^*=\left\{
\begin{array}{ll}
      0 & \text{ if } j=0\\
      x_{i,jt^2} & \text{ if } 1\leq j<m \\
      1 & \text{ if } j=m
\end{array} 
\right. 
\end{equation}

Now, let $S_{i,j}=S_i\cap \{(x,y)\in [0,1]^2:x_{i,j-1}^*<x\leq x_{i,j}^*\}$. That is, the $S_{i,j}$'s partition each strip into $m$ cells, each containing $t^2$ vertices of $G$. Now that we have our partition of the unit square, we can define our mapping of vertices in $G$ to vertices in {$L(m^2,d,t^2)$, where $d$ is yet to be determined.} The mapping is so that for each $(i,j)$, all points in $S_{i,j}$ are mapped injectively to the vertices of clique $C_{i,j}$ in any order.

Let $F:P\to [n]^2$, be a mapping from $P$ to the vertices of $L(m^2,d,t^2)$ as described above. By definition $F$ is a bijection. If $F$ is a homomorphism then {any orthogonal colouring of $L(m^2,d,t^2)$ induces an orthogonal colouring of $G$: assign each vertex $v$ of $G$ the colours of $F(v)$.} We can then apply Lemma \ref{Lemma: Orthogonal Colouring Lattice} to obtain an optimal orthogonal colouring of $L(m^2,d,t^2)$ and thus also of $G$. To {determine the correct value of $d$ and to } show that the mapping $F$ is a {homomorphism (with high probability)}, we need to bound the values of the variables $y_i^*$ and $x_{i,j}^*$. For this, we will use the following two lemmas, which follow directly from the Chernoff bound, stated earlier as Theorem \ref{ChernoffBound} .

If the points that determine the random geometric graph were evenly spaced in the unit square, then we would have that $y_i^*=\frac{i}{m}$ for all $i$. We can use the following lemma to establish that this is approximately true. In the following, we say that an event $E$ holds \emph{with exponential probability (w.e.p.)} if the probability that $E$ does not hold is exponentially small. Precisely, $E$ holds w.e.p.~if $P(E)=1-e^{-\Omega(-(\ln n)^2)}$.

\begin{lem}\label{Lemma: 1}
Let $P$ be a collection of $n=m^2t^2$ points, chosen u.a.r.~from the unit square.
If $\delta=\frac{2\ln n}{\sqrt{n}}=\frac{2\ln n}{mt}$, then with exponential probability, for all points $(a,b)\in P$ and for all integers $1\leq s<m$, if
$$\left|\{(x,y)\in P:y\leq b\}\right|= s mt^2,$$
then $$\left|b-\frac{s}{m}\right|\leq \delta.$$ 
\end{lem}

\begin{pf}
{
Fix $p_0=(a,b)\in P$ and $1\leq s<m$, and assume $b-\frac{s}{m}> \delta$. Let $X=|\{(x,y)\in P:y\leq b\}|$. Relative to point $(a,b)$ all other $n-1$ points are chosen independently u.a.r.~in the unit square. For all $p=(x,y)\in P\setminus \{ p_0\}$, let $X_p$  be the indicator variable of $\{ y\leq b\}$, and let $X=\sum_{p\in P\setminus\{ p_0\}} X_p$. Then ${\mathbb E}(X_p)=b$. Let $n$ be large enough so that $\delta n\geq 2$. Then
$$\mu={\mathbb E}(X)=b(n-1)>(\frac{s}{m}+\delta)n -1\geq smt^2 +\delta n/2=smt^2+mt\ln(n).$$
By the Chernoff bound, 
\[
{\mathbb P}(X=  smt^2)\leq {\mathbb P}(X<(1-\frac{mt\ln{n}}{\mu})\mu)\leq e^{-\frac{n(\ln{n})^2}{3\mu}}
\]
By definition, $\mu\leq n-1<n$, so the absolute value of the exponent is greater than $ \ln(n)^2/3$, and thus the probability that $X= smt^2$ is exponentially small. Thus w.e.p.~if $b-\frac{s}{m}> \delta$ then $\left|\{(x,y)\in P:y\leq b\}\right|\neq s mt^2$.  With a similar argument, we can show that, if $\frac{s}{m}-b> \delta$ then w.e.p.~we have that $\left|\{(x,y)\in P:y\leq b\}\right|\neq s mt^2$. 

For any $p=(a,b)\in P$, let $\mathcal{E}_{p}$ be the event that $\left|\{(x,y)\in P:y\leq b\}\right|= s mt^2$ and 
$ \left|b-\frac{s}{m}\right|> \delta$. We saw that ${\mathbb P}(\mathcal{E}_p)$ is exponentially small. Since $|P|=n$, we have by a union bound that ${\mathbb P}(\cup_{p\in P} \mathcal{E}_p)$ is also exponentially small. This concludes the argument. 
}
\end{pf}

Lemma \ref{Lemma: 1} says that w.e.p.~if the number of points with $y$-value less than or equal to $b$ is some multiple $s$ of $mt^2$, then $\left|b-\frac{s}{m}\right|\leq \delta =2\ln(n)/mt$. Since the $y_i^*$'s are the positions where this would occur, we obtain the following corollary.

\begin{cor}\label{Corollary: 1}
Consider $n=m^2t^2$ points chosen uniformly at random from the unit square, and let $y_i^*$, $0\leq i\leq m$ be defined as in \eqref{def: yi*}. With exponential probability, for all $1\leq i<m$, 
\begin{equation}
\label{yiBound}
|y_i^*-\frac{i}{m}|\leq \frac{2\ln n}{mt}.
\end{equation}
\end{cor}

Corollary \ref{Corollary: 1} gives an upper and lower bound on the position of $y_i^*$. This will be useful in determining the parameter $d$ and proving that the map $F$ is indeed a homomorphism.
%the proof for bounding the Euclidean distance between two points. 
We now proceed to bound the positions of the $x_{i,j}^*$.

\begin{lem}\label{Lemma: 2}{
Let $P$ be a collection of $n=m^2t^2$ points, chosen u.a.r.~from the unit square. %For all $p\in P$, let $\mathcal{E}_p$ be as defined in the proof of Lemma \ref{Lemma: 1}.
Let $\delta=\frac{2\ln n}{t^{1/2}n^{1/4}}=\frac{2\ln n}{t\sqrt{m}}$. Then with exponential probability, for all points $p_1=(a_1,b_1)$, $p_2=(a_2,b_2)$, $p_3=(a_3,b_3)\in P$, and for all integers $1\leq s<m$, 
\begin{eqnarray*}
%\left|\{(x,y)\in P:x\leq b_1,b_2<y\leq b_3\}\neq k^2\right|
&\text{if }&\left|\{(x,y)\in P:b_2<y\leq b_3\}\right|=t^2m\\
%&\text{and }&\left|\{(x,y)\in P:y\leq b_3\}\right|=(r+1)t^2m\\
&\text{and }&\left|\{(x,y)\in P:b_2<y\leq b_3,\,x\leq a_1\}\right|=st^2,\\
&\text{then }&
\left|a_1-\frac{s}{m}\right|\leq \delta.
\end{eqnarray*}
}
\end{lem}

\begin{pf}
{
Fix $p_i=(a_i,b_i)\in P$ for $i=1,2,3$ and fix $1\leq s<m$. Let $P_0=\{(x,y)\in P:b_2<y\leq b_3\}$ and assume $|P_0|=t^2m=tn$.
Also assume that $a_1-\frac{s}{m}> \delta$.

Let $X=|\{(x,y)\in P_0\setminus\{ p_1\}:x\leq a_1\}|$. Since $x$-values and $y$-values of the randomly chosen points are independent, $X$ is the sum of $t^2m-1$ independent Bernouilli variables with success probability $a_1$. Thus
$$
\mu={\mathbb E}(X)=a_1(t^2m-1)\geq a_1t^2m-1>(\frac{s}{m}+\delta)t^2m-1 \geq st^2+t\sqrt{m}\ln(n),
$$
for $n$ large enough so that $\delta t^2m\geq 2$.
By the Chernoff bound, 
\[
{\mathbb P}(X= st^2)\leq {\mathbb P}(X<(1-\frac{t\sqrt{m}\ln{n}}{\mu})\mu)\leq e^{-\frac{t^2m(\ln{n})^2}{3\mu}}
\]
By definition, $\mu\leq |P_0|=t^2m$, so the absolute value of the exponent is at least $ \ln(n)^2/3$, and thus the probability that $X= smt^2$ is exponentially small. Thus w.e.p.~if $b-\frac{s}{m}> \delta$ then $\left|\{(x,y)\in P:y\leq b\}\right|\neq s mt^2$.  With a similar argument, we can show that, if $\frac{s}{m}-b> \delta$ then w.e.p.~we have that $\left|\{(x,y)\in P:y\leq b\}\right|\neq s mt^2$. Thus, the implication holds for this particular choice of $p_1,p_2,p_3$ and $s$. The result then follows by a union bound.

}
\end{pf}

Lemma \ref{Lemma: 2} says that if $b_1$ and $b_3$ border a strip containing $t^2m$ points, and there are $st^2$ points in this strip with $x$-value less than $a_1$, then $|a_1-\frac{s}{m}|\leq \delta =2\ln{n}/t\sqrt{m}$. Since the $x_{i,j}^*$'s are the positions where this would occur, we obtain the following corollary.

\begin{cor}\label{Corollary: 2}
Consider $n^2=m^2t^2$ points chosen uniformly at random from the unit square, and let $x_{i,j}^*$, $0\leq i,j\leq m$ be defined as in \eqref{def:xij*}. With exponential probability, for all $1\leq i,j <m$, 
\begin{equation}\label{xijBound}
|x_{i,j}^*-\frac{j}{m}|\leq  \frac{2\ln n}{t\sqrt{m}}.
\end{equation}
\end{cor}

With these two results in hand, we now bound the Euclidean distance between any pair of points in $S_{i,j}$ and $S_{i',j' }$. 

{
\begin{cor}\label{Corollary: 3}
Assume that $y_i^*$ and $x_{i,j}^*$ satisfy bounds \eqref{yiBound} and \eqref{xijBound} for all $1\leq i,j\leq m$. Then for all points $p_1\in S_{i,j}$ and $p_2\in S_{i',j'}$ we have that
\begin{equation}\label{eq:pointbound}
\|p_1-p_2\|\geq \frac{1}{m}\left( \max\{ |i-i'|,|j-j'|\} - 1-\frac{4\sqrt{m}\ln{n}}{t}\right).
\end{equation}

\end{cor}

\begin{pf}
Let $p_1=(a_1,b_1) $ and $p_2=(a_2,b_2)$.
Each cell $S_{i,j}$ is bounded below and above by $y_{i-1}^*$, $y_i^*$, and left and right by
$x_{i,j-1}^*$ and $x^*_{i,j}$. Suppose first (by possible relabelling) that $b_2\geq b_1$, and thus $i'\geq i$. Then by Corollary \ref{Corollary: 1},
\[
\|p_1-p_2\|\geq b_2-b_1\geq  y^*_{i'-1}-y^*_i\geq (i'-1-i)/m - 2\left(\frac{2\ln{n}}{mt}\right).
%\geq \frac{|i'-i|-3}{m}.
\]
Since $\sqrt{m}\geq 1$, the bound \eqref{eq:pointbound} follows. 

Next assume \emph{wlog} that $a_1\leq a_2$ and thus $j'\geq j$ (where $i,i'$ can have any order). Then by Corollary \ref{Corollary: 2},
\[
\|p_1-p_2\|\geq |a_2-a_1|\geq | x^*_{i',j'-1}-x^*_{i,j}|\geq (j'-1-j)/m - 2\left(\frac{2\ln{n}}{t\sqrt{m}}\right).
%\geq \frac{|j'-j|-3}{m}.
\]
This matches the bound \eqref{eq:pointbound}.
\end{pf}
}
We can use this corollary directly to prove a lemma that establishes when the map $F$ is a homomorphism.

{
\begin{lem}
Let $n=m^2t^2$ where $1\leq t\leq m$, and assume $G\sim RG(n,r)$. If $d> rm   + 4\sqrt{m}\ln{n}/t$, then w.e.p.~there is a bijective homomorphism from $G$ into $L(m^2,d,t^2)$.
%As a consequence, if $mt$ is a perfect square and $r\leq cn^{-1/4} $ for some $c<1$, then a.a.s.~$G$ has an optimal orthogonal colouring. 
\end{lem}
\begin{pf}
Assume that $y_i^*$ and $x_{i,j}^*$ satisfy bounds \eqref{yiBound} and \eqref{xijBound} for all $1\leq i,j\leq m$. 
Let $F$ be the bijection from $G$ to $L(m^2,d,t)$ defined by mapping each set $S_{i,j}$ to the corresponding clique $C_{i,j}$ (in any order). It suffices to show that $F$ preserves edges. Note that for points $p_1,p_2\in P$,  $p_1$ is adjacent to $p_2$ in $G$ if and only if $\|p_1-p_2\|\leq r$. Therefore we need to prove that for all $p_1\in S_{i,j}$ and $p_2\in S_{i',j'}$, if $\|p_1-p_2\|\leq r$ then $\max\{ |i'-i|,|j'-j|\}\leq d$.  
We prove the converse. Assume that $\max\{ |i'-i|,|j'-j|\}\geq  d+1$. 
Then by Corollary \ref{Corollary: 3} we have that
\[
\| p_1-p_2\|\geq \frac{1}{m}\left(d-4\frac{\sqrt{m}\ln{n}}{t}\right)> r,
\]
where the last inequality follows from the condition on $d$. This shows that $F$ is a (bijective) homomorphism.
\end{pf}

The last step in establishing that $RG(n,r)$ has an optimal orthogonal colouring when $r<n^{-1/4}$ is to show that there are  parameters $m,d,t$  such that $L(m^2,d,t^2)$ has an optimal orthogonal colouring. 

\begin{thrm}
Fix $c\in (0,1)$ and assume $r\leq cn^{-1/4}$. Let $G\sim RG(n,r)$. If $n=m^2t^2$ where $\ln(n)^4\leq t\leq rn^{1/2}/\ln(n)$, 
then a.a.s.~$G$ has an optimal orthogonal colouring. 
\end{thrm}

\begin{pf}
Fix $c\in (0,1)$ and assume $r= r(n)= cn^{-1/4}$.
Take $d=\lfloor rm +4\sqrt{m}\ln(n)/t\rfloor+1$, so
\begin{equation}\label{eq:d_upper_bound}
d+1\leq rm\left(1+\frac{2}{rm}\right)+ \frac{4\sqrt{m}\ln n}{t}.
\end{equation}
Now since $mt=n^{1/2}$ and $t\leq rn^{1/2}/\ln{n}$,
\[
rm = rn^{1/2}/t\geq \ln(n)\rightarrow \infty\text{ as } n\rightarrow \infty,
%rt\sqrt{m}\geq c\ln{n}.
\]
and thus $2/rm=o(1)$. Fix $\epsilon >0$ and let $n$ be large enough so that 
\begin{equation}\label{eq:condition}
\left(1+\frac{2}{rm}\right)^2<\frac{1-\epsilon}{c^2},
\end{equation}
and thus
\[
\left(1+\frac{2}{rm}\right)^2 r^2m^2< (1-\epsilon)n^{-1/2}m^2.
\]
Then we have that
\begin{eqnarray*}
t(d+1)^2&\leq&t\left( rm\left(1+\frac{2}{rm}\right)+ \frac{4\sqrt{m}\ln n}{t}\right)^2\\
&<&(1-\epsilon)tn^{-1/2}m^2+ 8\sqrt{(1-\epsilon)}n^{-1/4}m^{3/2}\ln n+16m(\ln n)^2/t\\
&\leq & m(1-\epsilon+8t^{-1/2}+16/(\ln n)^2)\\
&=& m\left( 1-\epsilon +o\left(\frac{1}{\ln n}\right)\right),
\end{eqnarray*}
where the last two lines follow from the fact that $t\geq \ln(n)^4$.

Therefore, for large enough $n$ we have that $t(d+1)^2<m$. Then we are in Case 3 of Theorem \ref{Lemma: Orthogonal Colouring Lattice}, and $O\chi(G)=O\chi(L(m^2,d,t^2)=mt=\sqrt{n}$. Therefore, if  $y_i^*$ and $x_{i,j}^*$ satisfy bounds \eqref{yiBound} and \eqref{xijBound} then $G$ has an optimal orthogonal colouring. By Lemma \ref{Lemma: 1} and Lemma \ref{Lemma: 2}, this occurs w.e.p.
\end{pf}
}

{Note that Theorem \ref{Theorem: Random Main Result} only applies for values of $n$ that have a factor $t\geq (\ln n)^4$. There are, of course, infinitely many such values. We believe that the proof of the theorem can be adjusted to accommodate values of $n$ that do not have a large factor. In that case, we can divide the unit square in cells so that most cells contain $t^2$ points, but at most one square in each slice $S_i$ has less than $t^2$ points. We have not pursued this argument since it promised to be overly technical. }

\bibliographystyle{amsplain}

\bibliography{Orthogonal_Colourings_of_Random_Geometric_Graphs_References}
\end{document}